\documentclass[11pt,a4paper]{article}

\usepackage{inputenc}
\usepackage{amsmath}
\usepackage{bm}
\usepackage{bbold}
\usepackage{amsthm}
\usepackage{enumerate}

\usepackage{hyperref}
\usepackage{breakurl}

\title{Algebraic solution of minimax single-facility constrained location problems with Chebyshev and rectilinear distances}

\author{N. Krivulin\thanks{Faculty of Mathematics and Mechanics, Saint Petersburg State University, 28 Universitetsky Ave., St.~Petersburg, 198504, Russia, 
nkk@math.spbu.ru.}
\thanks{This work was supported in part by the Russian Foundation for Basic Research (grant No. 18-010-00723, 20-010-00145).}}

\date{}

\newtheorem{theorem}{Theorem}
\newtheorem{lemma}[theorem]{Lemma}
\newtheorem{corollary}[theorem]{Corollary}

\theoremstyle{definition}

\begin{document}

\maketitle

\begin{abstract}
We consider location problems to find the optimal sites of placement of a new facility, which minimize the maximum weighted Chebyshev or rectilinear distance to existing facilities under constraints on a feasible location domain. We examine Chebyshev location problems in multidimensional space to represent and solve the problems in the framework of tropical (idempotent) algebra, which deals with the theory and applications of semirings and semifields with idempotent addition. The solution approach involves formulating the problem as a tropical optimization problem, introducing a parameter that represents the minimum value of the objective function in the problem, and reducing the problem to a system of parametrized inequalities. The necessary and sufficient conditions for the existence of a solution to the system serve to evaluate the minimum, whereas all corresponding solutions of the system present a complete solution of the optimization problem. With this approach we obtain direct, exact solutions represented in a compact closed form which is appropriate for further analysis and straightforward computations with polynomial time complexity. The solutions of the Chebyshev problems are then used to solve location problems with rectilinear distance in the two-dimensional plane. The obtained solutions extend previous results on the Chebyshev and rectilinear location problems without weights and with less general constraints.
\\

\textbf{Key-Words:} tropical optimization, idempotent semifield, constrained optimization problem, single-facility location problem, Chebyshev and rectilinear distances
\\

\textbf{MSC (2020):} 90C24, 15A80, 90B85, 90C47, 90C48
\end{abstract}

\section{Introduction}

Location problems present an important research domain in optimization, which dates back to the XVII century and originates in the influential works of P.~Fermat, E.~Torricelli, J.~J.~Sylvester, J.~Steiner and A.~Weber. Many results achieved in this domain are recognized as notable contributions to various fields, such as operations research, computer science and engineering.

To solve location problems, which are formulated in different settings, a variety of analytical approaches and computational techniques exist, including methods of linear and mixed-integer linear programming, methods of discrete, combinatorial and graph optimization \cite{Sule2001Logistics,Klamroth2002Singlefacility,Farahani2009Facility,Eiselt2011Foundations,Laporte2015Location}. Another approach, which finds increasing application in solving some classes of optimization problems, is to use models and methods of tropical mathematics.  

Tropical (idempotent) mathematics deals with the theory and applications of semirings and semifields with idempotent addition (see, e.g., \cite{Golan2003Semirings,Heidergott2006Maxplus,Gondran2008Graphs,Mceneaney2010Maxplus}). It includes tropical optimization as a research area concerned with optimization problems that are formulated and solved in the framework of tropical mathematics. In many cases, tropical optimization problems can be solved directly in closed form under general assumptions, whereas other problems have only algorithmic solutions based on iterative numerical procedures. For a brief overview of tropical optimization problems, one can consult, e.g., \cite{Krivulin2014Tropical}.

As a solution framework, tropical mathematics is used in \cite{Cuninghamegreen1991Minimax,Cuninghamegreen1994Minimax} to handle one-dimensional minimax location problems on graphs. A similar algebraic approach based on the theory of max-separable functions is implemented in \cite{Zimmermann1991Minmax,Zimmermann1992Optimization,Hudec1993Aservice,Hudec1999Biobjective,Tharwat2010Oneclass} to solve constrained minimax location problems. Further examples include the solutions, given in terms of idempotent algebra in \cite{Krivulin2011Anextremal,Krivulin2014Complete,Krivulin2015Onanalgebraic,Krivulin2016Using,Krivulin2017Using}, to unconstrained and constrained minimax single-facility location problems with Chebyshev and rectilinear distances.  

In this paper, we consider location problems to find the optimal sites of placement of a new facility, which minimize the maximum weighted Chebyshev or rectilinear distance to existing facilities under constraints on a feasible location domain. For any two vectors $\bm{r}=(r_{1},\ldots,r_{n})^{T}$ and $\bm{s}=(s_{1},\ldots,s_{n})^{T}$ in the real space $\mathbb{R}^{n}$, the Chebyshev distance (maximum, dominance or $L_{\infty}$ metric) is given by
\begin{equation}
\mathrm{d}_{\infty}(\bm{r},\bm{s})
=
\max_{1\leq i\leq n}|r_{i}-s_{i}|
=
\max_{1\leq i\leq n}\max\{r_{i}-s_{i},s_{i}-r_{i}\}.
\label{E-dinfty}
\end{equation}

The rectilinear distance (Manhattan, rectangular, taxi-cab, city-block or $L_{1}$ metric) is defined as
\begin{equation}
\mathrm{d}_{1}(\bm{r},\bm{s})
=
\sum_{1\leq i\leq n}|r_{i}-s_{i}|
=
\sum_{1\leq i\leq n}\max\{r_{i}-s_{i},s_{i}-r_{i}\}.
\label{E-d1}
\end{equation}

Suppose that we are given $m$ points $\bm{p}_{i}=(p_{1i},\ldots,p_{ni})^{T}\in\mathbb{R}^{n}$, positive reals $w_{i}$ (weights) and $d_{i}$ (upper bounds), and reals $h_{i}$ (addends) for all $i=1,\ldots,m$. We need to locate a new point $\bm{x}=(x_{1},\ldots,x_{n})^{T}$ in a feasible location domain $S\subset\mathbb{R}^{n}$ to minimize the maximum weighted distance with addends, in the sense of a metric $\mathrm{d}$, from $\bm{x}$ to existing points under upper bound constraints on these distances. The problem is formulated in the form
\begin{equation}
\begin{aligned}
&
\text{min}
&&
\max_{1\leq j\leq m}(w_{j}\mathrm{d}(\bm{x},\bm{p}_{j})+h_{j});
\\
&
\text{s.~t.}
&&
\mathrm{d}(\bm{x},\bm{p}_{j})
\leq
d_{j},
\quad
j=1,\ldots,m;
\\
&&&
\bm{x}\in S.
\end{aligned}
\label{P-minwjdxrjhj-dxrjdj-xS}
\end{equation}

We examine the problem with Chebyshev and rectilinear distances under different settings of the dimension $n$ and of the feasible location area $S$. In the case of the Chebyshev distance, we retain the general setting of a real space of arbitrary dimension $n$. Given real numbers $b_{ij}$, $c_{i}\ne0$, $f_{i}$ and $g_{i}$ such that $f_{i}\leq g_{i}$ for all $i,j=1,\ldots,n$, the location area is described by the set
\begin{equation}
S
=
\{(x_{1},\ldots,x_{n})^{T}|\ b_{ik}+c_{k}x_{k}\leq c_{i}x_{i},\ f_{i}\leq x_{i}\leq g_{i},\ 1\leq i,k\leq n\}.
\label{E-S-bikxkckxicifixigi}
\end{equation}

This area takes the form of the intersection, if it exists, of the half-spaces defined by the inequalities $b_{ik}+c_{k}x_{k}\leq c_{i}x_{i}$ and of the hyper-rectangle given by the double inequalities $f_{i}\leq x_{i}\leq g_{i}$. We first solve the location problem under the condition that $c_{i}=1$ for all $i$ and then extend the solution to the problem for the case of arbitrary real $c_{i}\ne1$ for all $i=1,\ldots,n$.

We extend the set of real numbers by adding $-\infty$ and allow the parameters $b_{ij}$ to take the value $-\infty$. The setting $b_{ik}=-\infty$ formally makes the corresponding inequality constraint $b_{ik}+c_{k}x_{k}\leq c_{i}x_{i}$ inoperative and hence enables to define problems with reduced systems of constraints.

In the rectilinear case, we solve two more specific two-dimensional problems with different location areas defined on the real plane as follows. Given real numbers $f_{1}$, $f_{2}$, $g_{1}$, $g_{2}$, $a$ and $b$ such that $f_{1}\leq g_{1}$, $f_{2}\leq g_{2}$ and $a\leq b$, we first consider the location set
\begin{equation}
S
=
\{(x_{1},x_{2})^{T}|\ f_{1}-x_{2}\leq x_{1}\leq g_{1}-x_{2},\ f_{2}+x_{1}\leq x_{2}\leq g_{2}+x_{1},\ a\leq x_{1}\leq b\},
\label{E-S-f1x2x1g1x2f2x1x2g2x1ax1b}
\end{equation}
which presents the intersection of the $45^{\circ}$ tilted rectangle defined by the inequalities $f_{1}-x_{2}\leq x_{1}\leq g_{1}-x_{2}$ and $f_{2}+x_{1}\leq x_{2}\leq g_{2}+x_{1}$, and the vertical rectilinear strip area given by the boundary conditions $a\leq x_{1}\leq b$. 

Next, we introduce an additional real parameter $c\ne1$ and consider the feasible set 
\begin{equation}
S
=
\{(x_{1},x_{2})^{T}|\ f_{1}-x_{2}\leq x_{1}\leq g_{1}-x_{2},\ f_{2}+x_{1}\leq x_{2}\leq g_{2}+x_{1},\ a+x_{2}\leq cx_{1}\leq b+x_{2}\},
\label{E-S-f1x2x1g1x2f2x1x2g2x1ax2cx1bx2}
\end{equation}
which is the intersection of the $45^{\circ}$ tilted rectangle provided by the first two inequalities, and an arbitrary tilted strip given by the double inequality $a+x_{2}\leq cx_{1}\leq b+x_{2}$.     

Note that the above descriptions of the feasible sets in $n$-dimensional and two-dimensional spaces cover a range of geometrical configurations from nontrivial convex polytopes (polygons) to line segments.

The problems under consideration and their special cases are examined in many works which offer various solutions to the problem. First note that these problems can be formulated as linear programs and then solved using an appropriate linear programming computational procedure such as the simplex or Karmarkar algorithms. This approach, however, provides a numerical solution, if it exists, rather than a direct, complete solution in an exact analytical form.

For the unconstrained problems with rectilinear distance and equal weights, direct explicit solutions are obtained in \cite{Elzinga1972Geometrical,Francis1972Ageometrical} using geometric arguments. A solution for the weighted problem with rectilinear distance is given in \cite{Dearing1972Onsome}, which involves decomposition into independent one-dimensional subproblems solved by a reduction to equivalent network flow problems. In \cite{Krivulin2011Anextremal,Krivulin2014Complete,Krivulin2015Onanalgebraic,Krivulin2016Using,Krivulin2017Using}, an approach based on idempotent algebra is applied to solve unweighted unconstrained and constrained location problems. Further results on both unweighted and weighted location problems can be found in the survey papers \cite{Francis1983Locational,Brandeau1989Anoverview,ReVelle2005Location,Brimberg2009Optimizing,Chhajed2013Facility}, as well as in the books \cite{Sule2001Logistics,Klamroth2002Singlefacility,Farahani2009Facility,Eiselt2011Foundations,Laporte2015Location}. 

In this paper, we represent and examine the location problems in the framework of tropical (idempotent) algebra. We start with the solution of location problems with Chebyshev distance in multidimensional space. The solution approach follows the analytical technique developed in \cite{Krivulin2014Complete,Krivulin2015Multidimensional,Krivulin2016Using,Krivulin2017Using}, which involves formulating the problem as a tropical optimization problem, introducing a parameter that represents the minimum value in the problem, and reducing the problem to a system of parametrized inequalities. The necessary and sufficient conditions for the existence of a solution of the system serve to evaluate the minimum, whereas all corresponding solutions of the system present a complete solution of the optimization problem. With this approach, we obtain direct, exact solutions represented in a compact closed form which is appropriate for further analysis and straightforward computations with polynomial time complexity. The solutions of the Chebyshev problems are then used to solve location problems with rectilinear distance in the two-dimensional plane.

The obtained solutions further extend previous work in \cite{Krivulin2014Complete,Krivulin2017Using} on the location problems without weights (positive equal-weighted problems). Furthermore, the presented results enlarge the findings published in the short conference paper \cite{Krivulin2018Algebraic} by expanding the class of feasible location sets under consideration. The new solutions which are given in an explicit form can serve to supplement and complement existing methods and be of particular interest when the application of known numerical algorithmic solutions, for one reason or other, appears to be impractical or impossible.

\section{Elements of Tropical Algebra}

In this section, we present a brief introduction to tropical (idempotent) algebra to provide a formal analytical framework for the solution of the location problems in the sequel. For further details on the theory and applications of tropical mathematics, one can refer, for example, to recent works \cite{Golan2003Semirings,Heidergott2006Maxplus,Gondran2008Graphs,Mceneaney2010Maxplus}.

An idempotent semifield is an algebraic system $(\mathbb{X},\oplus,\otimes,\mathbb{0},\mathbb{1})$ where $\mathbb{X}$ is a nonempty set that has distinct elements $\mathbb{0}$ (zero) and $\mathbb{1}$ (one). The set $\mathbb{X}$ is equipped with binary operations $\oplus$ (addition) and $\otimes$ (multiplication) such that $(\mathbb{X},\oplus,\mathbb{0})$ is a commutative idempotent monoid, $(\mathbb{X}\setminus\{\mathbb{0}\},\otimes,\mathbb{1})$ is an Abelian group, and $\otimes$ distributes over $\oplus$.

In an idempotent semifield, addition is idempotent, which means that $x\oplus x=x$ for all $x\in\mathbb{X}$, and induces a partial order by the rule: $x\leq y$ if and only if $x\oplus y=y$. This order is assumed to constitute a total order on $\mathbb{X}$. With respect to this order, the operations $\oplus$ and $\otimes$ are monotone, which implies that the inequality $x\leq y$ results in $x\oplus z\leq y\oplus z$ and $x\otimes z\leq y\otimes z$. Furthermore, addition has an extremal property in the sense that the inequalities $x\leq x\oplus y$ and $y\leq x\oplus y$ hold for all $x,y\in\mathbb{X}$. Finally, the inequality $x\oplus y\leq z$ is equivalent to the system of inequalities $x\leq z$ and $y\leq z$.

Multiplication is invertible, which provides for each $x\ne\mathbb{0}$ a unique inverse $x^{-1}$ such that $x\otimes x^{-1}=\mathbb{1}$. Inversion is antitone to turn the inequality $x\leq y$ where $x,y\ne\mathbb{0}$ into $x^{-1}\geq y^{-1}$. In what follows, the multiplication sign $\otimes$ is, as usual, omitted to save writing.

Integer powers are used in the standard way to indicate iterated products: $x^{0}=\mathbb{1}$, $x^{p}=x\otimes x^{p-1}$, $x^{-p}=(x^{-1})^{p}$ and $\mathbb{0}^{p}=\mathbb{0}$ for all $x\in\mathbb{X}$ and integer $p>0$. Furthermore, the equation $x^{p}=a$ has the unique solution $x=a^{1/p}$ for each $a\in\mathbb{X}$ and integer $p>0$, which allows the powers to have rational exponents. Moreover, it is assumed that the power notation can be further extended to real exponents (e.g., by the usual extra limiting process) to have real powers and the power rules well defined. Exponentiation is monotone, which means that the inequality $x\leq y$ yields $x^{p}\leq y^{p}$ if $p>0$ and $x^{p}\geq y^{p}$ if $p<0$.   

An analogue of the binomial identity holds in the form $(a\oplus b)^{r}=a^{r}\oplus b^{r}$ for any $a,b\in\mathbb{X}$ and nonnegative real $r$.

As an example, we consider the real semifield $(\mathbb{R}\cup\{-\infty\},\max,+,-\infty,0)$, also known as the $(\max,+)$-algebra, where $\oplus=\max$, $\otimes=+$, $\mathbb{0}=-\infty$ and $\mathbb{1}=0$. In the $(\max,+)$-algebra, the power $x^{y}$ coincides with the usual arithmetic product $xy$, and the inverse $x^{-1}$ with the opposite number $-x$. The order induced by the idempotent addition corresponds to the natural linear order on $\mathbb{R}$. This semifield is used below for modeling the location problems under investigation. 

The algebra of vectors and matrices over idempotent semifields is introduced in the ordinary way. The vector (matrix) operations follow the conventional rules where the operations $\oplus$ and $\otimes$ are used instead of arithmetic addition and multiplication. In the following, all vectors are considered column vectors unless transposed to row vectors. A vector that has all elements equal to $\mathbb{0}$ is the zero vector. A vector with all elements not equal to $\mathbb{0}$ is called regular. 

A square matrix that has all entries equal to $\mathbb{1}$ on the diagonal and to $\mathbb{0}$ everywhere else, is the identity matrix denoted by $\bm{I}$. For any square matrix $\bm{A}$ and positive integer $p$, the power notation indicates iterated matrix products $\bm{A}^{0}=\bm{I}$, $\bm{A}^{p}=\bm{A}\bm{A}^{p-1}$. For any $(n\times n)$-matrix $\bm{A}=(a_{ij})$, the trace is given by
\begin{equation*}
\mathop\mathrm{tr}\bm{A}
=
a_{11}\oplus\cdots\oplus a_{nn}.
\end{equation*}

The properties of the scalar operations $\oplus$ and $\otimes$ with respect to the order relation $\leq$ are readily extended to the vector (matrix) operations where the inequalities are considered componentwise. For any nonzero vector $\bm{x}=(x_{i})$, the multiplicative conjugate transpose is a row vector $\bm{x}^{-}=(x_{i}^{-})$ with the elements $x_{i}^{-}=x_{i}^{-1}$ if $x_{i}\ne\mathbb{0}$, and $x_{i}^{-}=\mathbb{0}$ otherwise. For any regular vectors $\bm{x}$ and $\bm{y}$ such that $\bm{x}\leq\bm{y}$, the conjugate transposition yields $\bm{x}^{-}\geq\bm{y}^{-}$.

We conclude the overview with two results of tropical linear algebra. First suppose that, given an $(m\times n)$-matrix $\bm{A}$ and an $m$-vector $\bm{d}$, we need to find all $n$-vectors $\bm{x}$ that are solutions of the inequality 
\begin{equation}
\bm{A}\bm{x}
\leq
\bm{d}.
\label{I-Axleqd}
\end{equation}

A complete solution of inequality \eqref{I-Axleqd} is given as follows (see, e.g., \cite{Cuninghamegreen1994Minimax,Heidergott2006Maxplus}).
\begin{lemma}\label{L-Axleqd}
Let $\bm{A}$ be a matrix without zero columns and $\bm{d}$ a regular vector. Then, all solutions of inequality \eqref{I-Axleqd} are given by $\bm{x}\leq(\bm{d}^{-}\bm{A})^{-}$.
\end{lemma}

Furthermore, given an $(n\times n)$-matrix $\bm{A}$ and an $n$-vector $\bm{b}$, we consider the problem to find all regular $n$-vectors $\bm{x}$ that satisfy the inequality
\begin{equation}
\bm{A}\bm{x}
\oplus
\bm{b}
\leq
\bm{x}.
\label{I-Axbleqx}
\end{equation}

To describe a solution to the problem, we introduce a function that maps any $(n\times n)$-matrix $\bm{A}$ onto the scalar
\begin{equation*}
\mathop\mathrm{Tr}(\bm{A})
=
\mathop\mathrm{tr}\bm{A}
\oplus\cdots\oplus
\mathop\mathrm{tr}\bm{A}^{n}.
\end{equation*}

Provided that $\mathop\mathrm{Tr}(\bm{A})\leq\mathbb{1}$, the asterate operator (the Kleene star) transforms the matrix $\bm{A}$ into the matrix
\begin{equation*}
\bm{A}^{\ast}
=
\bm{I}\oplus\bm{A}\oplus\cdots\oplus\bm{A}^{n-1}.
\end{equation*}

The next statement presents a solution proposed in \cite{Krivulin2015Multidimensional} to inequality \eqref{I-Axbleqx}.
\begin{theorem}\label{T-Axbleqx}
For any matrix $\bm{A}$ and vector $\bm{b}$, the following statements hold:
\begin{enumerate}
\item
If $\mathop\mathrm{Tr}(\bm{A})\leq\mathbb{1}$, then all regular solutions of inequality \eqref{I-Axbleqx} are given by $\bm{x}=\bm{A}^{\ast}\bm{u}$ where the parameter vector $\bm{u}$ satisfies the condition $\bm{u}\geq\bm{b}$.
\item
If $\mathop\mathrm{Tr}(\bm{A})>\mathbb{1}$, then there are no regular solutions.
\end{enumerate}
\end{theorem}

Below, we represent the location problems under study in terms of idempotent algebra and obtain direct, complete solutions of the problems.

\section{Location with Chebyshev Distance}

We start with a solution of the location problem defined on the $n$-dimensional vector space with Chebyshev metric \eqref{E-dinfty}. In the framework of $(\max,+)$-algebra, the Chebyshev distance between vectors $\bm{r}=(r_{i})$ and $\bm{s}=(s_{i})$ in $\mathbb{R}^{n}$ is given by
\begin{equation*}
\mathrm{d}_{\infty}(\bm{r},\bm{s})
=
\bigoplus_{1\leq i\leq n}(s_{i}^{-1}r_{i}\oplus r_{i}^{-1}s_{i})
=
\bm{s}^{-}\bm{r}\oplus\bm{r}^{-}\bm{s}.
\end{equation*}

The objective function in problem \eqref{P-minwjdxrjhj-dxrjdj-xS} takes the form
\begin{equation*}
\bigoplus_{1\leq j\leq m}h_{j}(\bm{p}_{j}^{-}\bm{x}\oplus\bm{x}^{-}\bm{p}_{j})^{w_{j}},
\end{equation*}
whereas the upper bound constraints are written as
\begin{equation*}
\bm{p}_{j}^{-}\bm{x}
\oplus
\bm{x}^{-}\bm{p}_{j}
\leq
d_{j},
\quad
j=1,\ldots,m.
\end{equation*}

Furthermore, the feasible location area which is defined by \eqref{E-S-bikxkckxicifixigi} becomes
\begin{equation}
S
=
\{(x_{1},\ldots,x_{n})^{T}|\ b_{ik}x_{k}^{c_{k}}\leq x_{i}^{c_{i}},\ f_{i}\leq x_{i}\leq g_{i},\ 1\leq i,k\leq n\}.
\label{E-S-bikxkckxici_1}
\end{equation}

To solve the problem, we first examine a particular case where the feasible set $S$ is defined under the less general assumption that $c_{i}=1$ for all $i=1,\ldots,n$. The solution for this case demonstrates the key points of the approach and yields results in a more compact vector form. The solution procedure is then exploited as a template to handle the problem with arbitrary $c_{i}\ne1$.

\subsection{Solution for the Particular Case}

We start with a problem with a location set that is derived from \eqref{E-S-bikxkckxici_1} by setting $c_{i}=1$ for all $i=1,\ldots,n$, and given by 
\begin{equation*}
S
=
\{(x_{1},\ldots,x_{n})^{T}|\ b_{ik}x_{k}\leq x_{i},\ f_{i}\leq x_{i}\leq g_{i},\ 1\leq i,k\leq n\}.
\end{equation*}

To represent the constraints of the set $S$ in a compact vector form, we first combine the inequalities $b_{ik}x_{k}\leq x_{i}$ for all $k=1,\ldots,n$ into one equivalent inequality 
\begin{equation*}
b_{i1}x_{1}\oplus\cdots\oplus b_{in}x_{n}
\leq
x_{i},
\qquad
i=1,\ldots,n.
\end{equation*}

With the matrix and vector notation
\begin{equation*}
\bm{B}
=
(b_{ik}),
\qquad
\bm{f}
=
(f_{i}),
\qquad
\bm{g}
=
(g_{i}),
\end{equation*}
we describe the location area by the system of vector inequalities
\begin{equation*}
\bm{B}\bm{x}
\leq
\bm{x},
\qquad
\bm{f}
\leq
\bm{x}
\leq\bm{g}.
\end{equation*}

After substitution of the Chebyshev metric and vector description of the location area in terms of $(\max,+)$-algebra into problem \eqref{P-minwjdxrjhj-dxrjdj-xS}, we formulate the problem as follows:
\begin{equation}
\begin{aligned}
&
\text{min}
&&
\bigoplus_{1\leq j\leq m}h_{j}(\bm{p}_{j}^{-}\bm{x}\oplus\bm{x}^{-}\bm{p}_{j})^{w_{j}};
\\
&
\text{s.~t.}
&&
\bm{p}_{j}^{-}\bm{x}
\oplus
\bm{x}^{-}\bm{p}_{j}
\leq
d_{j},
\quad
j=1,\ldots,m;
\\
&&&
\bm{B}\bm{x}
\leq
\bm{x},
\quad
\bm{f}
\leq
\bm{x}
\leq
\bm{g}.
\end{aligned}
\label{P-minhjrjxxrjwj-rjxxrjdj-Bxx-fxg}
\end{equation}

Note that we assume all data involved in the formulation of problem~\eqref{P-minwjdxrjhj-dxrjdj-xS} to be real numbers, with the exception of the entries of the matrix $\bm{B}$, which can be set to the tropical zero $\mathbb{0}=-\infty$. Specifically, both the known vectors $\bm{p}_{j}$ for all $j=1,\ldots,m$, and the unknown vector $\bm{x}$ are considered regular.

To solve the problem obtained, we first introduce an additional parameter to represent the minimum value of the objective function, and then reduce the problem to a parametrized system of inequalities. Subsequently, we use existence conditions for solutions of the system to evaluate the value of the parameter. Finally, all solutions of the system, which correspond to this value, serve as a complete solution to the initial optimization problem. 

Let us denote the minimum value of the objective function by $\theta$. Then, all solutions of the problem must satisfy the equation
\begin{equation*}
\bigoplus_{1\leq j\leq m}h_{j}(\bm{p}_{j}^{-}\bm{x}\oplus\bm{x}^{-}\bm{p}_{j})^{w_{j}}
=
\theta.
\end{equation*}

Since we assume $\theta$ to be the minimum of the objective function, the set of solutions remains unchanged after replacing the equation by the inequality 
\begin{equation*}
\bigoplus_{1\leq j\leq m}h_{j}(\bm{p}_{j}^{-}\bm{x}\oplus\bm{x}^{-}\bm{p}_{j})^{w_{j}}
\leq
\theta.
\end{equation*}

Using the extremal property of idempotent addition, we replace the last inequality by an equivalent system of inequalities, which describes all solutions of problem \eqref{P-minhjrjxxrjwj-rjxxrjdj-Bxx-fxg} as follows:
\begin{equation}
\begin{aligned}
h_{j}(\bm{p}_{j}^{-}\bm{x}\oplus\bm{x}^{-}\bm{p}_{j})^{w_{j}}
&\leq
\theta,
\quad
j=1,\ldots,m;
\\
\bm{p}_{j}^{-}\bm{x}
\oplus
\bm{x}^{-}\bm{p}_{j}
&\leq
d_{j},
\quad
j=1,\ldots,m;
\\
\bm{B}\bm{x}
&\leq
\bm{x},
\\
\bm{f}
&\leq
\bm{x}
\leq
\bm{g}.
\end{aligned}
\label{S-hjrjxxrjwjtheta-rjxxrjdj-Bxx-fxg} 
\end{equation}

We use the tropical analogue of the binomial identity to replace the inequality $h_{j}(\bm{p}_{j}^{-}\bm{x}\oplus\bm{x}^{-}\bm{p}_{j})^{w_{j}}\leq\theta$ by the inequalities $h_{j}(\bm{p}_{j}^{-}\bm{x})^{w_{j}}\leq\theta$ and $h_{j}(\bm{x}^{-}\bm{p}_{j})^{w_{j}}\leq\theta$. Since exponentiation is monotone, these inequalities can be further rewritten by the usual power rules as $h_{j}^{1/w_{j}}\bm{p}_{j}^{-}\bm{x}\leq\theta^{1/w_{j}}$ and $h_{j}^{1/w_{j}}\bm{x}^{-}\bm{p}_{j}\leq\theta^{1/w_{j}}$, and then represented as the inequalities $\bm{p}_{j}^{-}\bm{x}\leq\theta^{1/w_{j}}h_{j}^{-1/w_{j}}$ and $\bm{x}^{-}\bm{p}_{j}\leq\theta^{1/w_{j}}h_{j}^{-1/w_{j}}$. Application of Lemma~\ref{L-Axleqd} to solve the first inequality with respect to $\bm{x}$ yields $\bm{x}\leq\theta^{1/w_{j}}h_{j}^{-1/w_{j}}\bm{p}_{j}$.

We again use Lemma~\ref{L-Axleqd} to solve the second inequality with respect to $\bm{p}_{j}$ and then multiply both sides of the result by $\theta^{-1/w_{j}}h_{j}^{1/w_{j}}$ to obtain the inequality $\theta^{-1/w_{j}}h_{j}^{1/w_{j}}\bm{p}_{j}\leq\bm{x}$. Finally, we combine the results into the double inequality $\theta^{-1/w_{j}}h_{j}^{1/w_{j}}\bm{p}_{j}\leq\bm{x}\leq\theta^{1/w_{j}}h_{j}^{-1/w_{j}}\bm{p}_{j}$.

In the same way, we replace the inequality $\bm{p}_{j}^{-}\bm{x}\oplus\bm{x}^{-}\bm{p}_{j}\leq d_{j}$ by the inequalities $\bm{x}\leq d_{j}\bm{p}_{j}$ and $d_{j}^{-1}\bm{p}_{j}\leq\bm{x}$ and represent them as $d_{j}^{-1}\bm{p}_{j}\leq\bm{x}\leq d_{j}\bm{p}_{j}$. 

We now rewrite system \eqref{S-hjrjxxrjwjtheta-rjxxrjdj-Bxx-fxg} in the form
\begin{equation*}
\begin{aligned}
\theta^{-1/w_{j}}h_{j}^{1/w_{j}}\bm{p}_{j}
&\leq
\bm{x}
\leq
\theta^{1/w_{j}}h_{j}^{-1/w_{j}}\bm{p}_{j},
\quad
j=1,\ldots,m;
\\
d_{j}^{-1}\bm{p}_{j}
&\leq
\bm{x}
\leq
d_{j}\bm{p}_{j},
\quad
j=1,\ldots,m;
\\
\bm{B}\bm{x}
&\leq
\bm{x},
\\
\bm{f}
&\leq
\bm{x}
\leq
\bm{g}.
\end{aligned}
\end{equation*}

Furthermore, we combine the left inequalities for all $j=1,\ldots,m$ into one, which provides a lower bound for $\bm{x}$. Next, we represent the right inequalities as $\bm{x}^{-}\geq\theta^{-1/w_{j}}h_{j}^{1/w_{j}}\bm{p}_{j}^{-}$, $\bm{x}^{-}\geq d_{j}^{-1}\bm{p}_{j}^{-}$, and $\bm{x}^{-}\geq\bm{g}^{-}$. Summing up these inequalities for all $j$ and conjugate-transposing the result yield an upper bound.

Thus, we have the double inequality
\begin{multline*}
\bm{B}\bm{x}
\oplus
\bigoplus_{1\leq j\leq m}
\theta^{-1/w_{j}}h_{j}^{1/w_{j}}
\bm{p}_{j}
\oplus
\bigoplus_{1\leq j\leq m}
d_{j}^{-1}
\bm{p}_{j}
\oplus
\bm{f}
\leq
\bm{x}
\\
\leq
\left(
\bigoplus_{1\leq j\leq m}
\theta^{-1/w_{j}}h_{j}^{1/w_{j}}
\bm{p}_{j}^{-}
\oplus
\bigoplus_{1\leq j\leq m}
d_{j}^{-1}
\bm{p}_{j}^{-}
\oplus
\bm{g}^{-}
\right)^{-}.
\end{multline*}

To simplify further formulas, we introduce the notation
\begin{equation}
\begin{aligned}
\bm{q}
&=
\bigoplus_{1\leq j\leq m}
\theta^{-1/w_{j}}h_{j}^{1/w_{j}}
\bm{p}_{j},
&
\bm{r}^{-}
&=
\bigoplus_{1\leq j\leq m}
\theta^{-1/w_{j}}h_{j}^{1/w_{j}}
\bm{p}_{j}^{-},
\\
\bm{s}
&=
\bigoplus_{1\leq j\leq m}
d_{j}^{-1}
\bm{p}_{j}
\oplus
\bm{f},
&
\bm{t}^{-}
&=
\bigoplus_{1\leq j\leq m}
d_{j}^{-1}
\bm{p}_{j}^{-}
\oplus
\bm{g}^{-}.
\end{aligned}
\label{E-p-q-f-g}
\end{equation}

In terms of this notation, the double inequality becomes
\begin{equation}
\bm{B}\bm{x}
\oplus
(\bm{q}\oplus\bm{s})
\leq
\bm{x}
\leq
(\bm{r}^{-}\oplus\bm{t}^{-})^{-}.
\label{I-Bxps-x-xqt}
\end{equation}

It follows from Theorem~\ref{T-Axbleqx} that the left inequality has regular solutions $\bm{x}$ if the condition $\mathop\mathrm{Tr}(\bm{B})\leq\mathbb{1}$ is valid, and admits only a trivial solution $\bm{x}=\bm{0}$ otherwise. Assuming that this condition holds, we apply Theorem~\ref{T-Axbleqx} to solve the left inequality in the parametric form $\bm{x}=\bm{B}^{\ast}\bm{u}$, where the vector of parameters $\bm{u}$ satisfies the condition $\bm{u}\geq\bm{q}\oplus\bm{s}$.

Next, we substitute $\bm{x}$ by $\bm{B}^{\ast}\bm{u}$ in the right inequality and apply Lemma~\ref{L-Axleqd} to solve this inequality for $\bm{u}$ in the form $\bm{u}\leq((\bm{r}^{-}\oplus\bm{t}^{-})\bm{B}^{\ast})^{-}$. As a result, we represent the solution to \eqref{I-Bxps-x-xqt} as
\begin{equation*}
\bm{x}
=
\bm{B}^{\ast}\bm{u},
\qquad
\bm{q}\oplus\bm{s}
\leq
\bm{u}
\leq
((\bm{r}^{-}\oplus\bm{t}^{-})\bm{B}^{\ast})^{-}.
\end{equation*}

The set of parameter vectors $\bm{u}$ defined by the obtained lower and upper bounds is nonempty if and only if the following condition holds:
\begin{equation*}
\bm{q}\oplus\bm{s}
\leq
((\bm{r}^{-}\oplus\bm{t}^{-})\bm{B}^{\ast})^{-}.
\end{equation*}

We now use the last inequality to evaluate the parameter $\theta$ and to refine the consistency condition for the constraints. We note that the row-vector $(\bm{r}^{-}\oplus\bm{t}^{-})\bm{B}^{\ast}$ is regular since both vectors $\bm{r}$ and $\bm{t}$ are regular, and $\bm{B}^{\ast}\geq\bm{I}$. Then, it follows from Lemma~\ref{L-Axleqd} that the inequality under consideration is the solution, with respect to the vector $\bm{q}\oplus\bm{s}$, of the inequality
\begin{equation*}
(\bm{r}^{-}\oplus\bm{t}^{-})\bm{B}^{\ast}(\bm{q}\oplus\bm{s})
\leq
\mathbb{1},
\end{equation*}
which we further rewrite as the system of inequalities
\begin{equation*}
\bm{r}^{-}\bm{B}^{\ast}\bm{q}
\leq
\mathbb{1},
\qquad
\bm{r}^{-}\bm{B}^{\ast}\bm{s}
\leq
\mathbb{1},
\qquad
\bm{t}^{-}\bm{B}^{\ast}\bm{q}
\leq
\mathbb{1},
\qquad
\bm{t}^{-}\bm{B}^{\ast}\bm{s}
\leq
\mathbb{1}.
\end{equation*}

We expand the first three inequalities into the inequalities
\begin{align*}
\left(
\bigoplus_{1\leq j\leq m}
\theta^{-1/w_{j}}h_{j}^{1/w_{j}}
\bm{p}_{j}^{-}
\right)
\bm{B}^{\ast}
\left(
\bigoplus_{1\leq l\leq m}
\theta^{-1/w_{l}}h_{l}^{1/w_{l}}
\bm{p}_{l}
\right)
&\leq
\mathbb{1},
\\
\left(
\bigoplus_{1\leq j\leq m}
\theta^{-1/w_{j}}h_{j}^{1/w_{j}}
\bm{p}_{j}^{-}
\right)
\bm{B}^{\ast}
\bm{s}
&\leq
\mathbb{1},
\\
\bm{t}^{-}
\bm{B}^{\ast}
\left(
\bigoplus_{1\leq l\leq m}
\theta^{-1/w_{l}}h_{l}^{1/w_{l}}
\bm{p}_{l}
\right)
&\leq
\mathbb{1},
\\
\bm{t}^{-}\bm{B}^{\ast}\bm{s}
&\leq
\mathbb{1},
\end{align*}
which are broken down into the inequalities
\begin{align*}
(
\theta^{-1/w_{j}}h_{j}^{1/w_{j}}
\bm{p}_{j}^{-}
)
\bm{B}^{\ast}
(
\theta^{-1/w_{l}}h_{l}^{1/w_{l}}
\bm{p}_{l}
)
&\leq
\mathbb{1},
\quad
j,l=1,\ldots,m;
\\
(
\theta^{-1/w_{j}}h_{j}^{1/w_{j}}
\bm{p}_{j}^{-}
)
\bm{B}^{\ast}
\bm{s}
&\leq
\mathbb{1},
\quad
j=1,\ldots,m;
\\
\bm{t}^{-}
\bm{B}^{\ast}
(
\theta^{-1/w_{l}}h_{l}^{1/w_{l}}
\bm{p}_{l}
)
&\leq
\mathbb{1},
\quad
l=1,\ldots,m;
\\
\bm{t}^{-}\bm{B}^{\ast}\bm{s}
&\leq
\mathbb{1}.
\end{align*}

Solving the first three inequalities with respect to $\theta$ yields the result
\begin{align*}
h_{j}^{\frac{w_{l}}{w_{j}+w_{l}}}
h_{l}^{\frac{w_{j}}{w_{j}+w_{l}}}
(
\bm{p}_{j}^{-}
\bm{B}^{\ast}
\bm{p}_{l}
)^{\frac{w_{j}w_{l}}{w_{j}+w_{l}}}
&\leq
\theta,
\quad
j,l=1,\ldots,m;
\\
h_{j}
(
\bm{p}_{j}^{-}
\bm{B}^{\ast}
\bm{s}
)^{w_{j}}
&\leq
\theta,
\quad
j=1,\ldots,m;
\\
h_{l}
(
\bm{t}^{-}
\bm{B}^{\ast}
\bm{p}_{l}
)^{w_{l}}
&\leq
\theta,
\quad
l=1,\ldots,m;
\\
\bm{t}^{-}\bm{B}^{\ast}\bm{s}
&\leq
\mathbb{1}.
\end{align*}

By combining the first three inequalities for each $j,l=1,\ldots,m$, we obtain the system
\begin{gather*}
\theta
\geq
\bigoplus_{1\leq j,l\leq m}
\left(
h_{j}^{\frac{w_{l}}{w_{j}+w_{l}}}
h_{l}^{\frac{w_{j}}{w_{j}+w_{l}}}
(
\bm{p}_{j}^{-}
\bm{B}^{\ast}
\bm{p}_{l}
)^{\frac{w_{j}w_{l}}{w_{j}+w_{l}}}
\oplus
h_{j}
(
\bm{p}_{j}^{-}
\bm{B}^{\ast}
\bm{s}
)^{w_{j}}
\oplus
h_{l}
(
\bm{t}^{-}
\bm{B}^{\ast}
\bm{p}_{l}
)^{w_{l}}
\right),
\\
\bm{t}^{-}\bm{B}^{\ast}\bm{s}
\leq
\mathbb{1}.
\end{gather*}

Consider the first inequality, which gives a lower bound for the parameter $\theta$. Since $\theta$ is assumed to represent the minimum in the problem, we set it equal to the right-hand side of the inequality.

As one can see, under the condition $\mathop\mathrm{Tr}(\bm{B})\leq\mathbb{1}$, the second inequality implies that the set of vectors $\bm{x}$ satisfying the double inequality $\bm{s}\leq\bm{x}\leq\bm{t}$ has a nonempty intersection with the set of $\bm{x}$ given by $\bm{B}\bm{x}\leq\bm{x}$. Therefore, both inequalities $\mathop\mathrm{Tr}(\bm{B})\leq\mathbb{1}$ and $\bm{t}^{-}\bm{B}^{\ast}\bm{s}\leq\mathbb{1}$ serve as necessary and sufficient conditions for the constraints of the problem to be consistent.

We now summarize the results obtained in the following statement in terms of an arbitrary idempotent semifield.
\begin{theorem}\label{T-minhjrjxxrjwj-rjxxrjdj-Bxx-fxg}
With notation \eqref{E-p-q-f-g}, suppose that the following conditions hold:
\begin{enumerate}
\item
$\mathop\mathrm{Tr}(\bm{B})\leq\mathbb{1}$,
\item
$\bm{t}^{-}\bm{B}^{\ast}\bm{s}\leq\mathbb{1}$.
\end{enumerate}

Then, the minimum value in problem \eqref{P-minhjrjxxrjwj-rjxxrjdj-Bxx-fxg} is equal to
\begin{equation}
\theta
=
\bigoplus_{1\leq j,l\leq m}
\left(
h_{j}^{\frac{w_{l}}{w_{j}+w_{l}}}
h_{l}^{\frac{w_{j}}{w_{j}+w_{l}}}
(
\bm{p}_{j}^{-}
\bm{B}^{\ast}
\bm{p}_{l}
)^{\frac{w_{j}w_{l}}{w_{j}+w_{l}}}
\oplus
h_{j}
(
\bm{p}_{j}^{-}
\bm{B}^{\ast}
\bm{s}
)^{w_{j}}
\oplus
h_{l}
(
\bm{t}^{-}
\bm{B}^{\ast}
\bm{p}_{l}
)^{w_{l}}
\right).
\label{E-theta}
\end{equation}

All solutions are given in parametric form by
\begin{equation*}
\bm{x}
=
\bm{B}^{\ast}\bm{u},
\end{equation*}
where the vector of parameters $\bm{u}$ satisfies the condition
\begin{equation*} 
\bm{q}\oplus\bm{s}
\leq
\bm{u}
\leq
((\bm{r}^{-}\oplus\bm{t}^{-})\bm{B}^{\ast})^{-}.
\end{equation*}
\end{theorem}

It is not difficult to see that the solution given by Theorem~\ref{T-minhjrjxxrjwj-rjxxrjdj-Bxx-fxg} can be computed in polynomial time in the number of points $m$ and the dimension of space $n$. Clearly, the computationally most demanding part of the solution is the calculation of the parameter $\theta$ according to \eqref{E-theta}. The evaluation of $\theta$ requires calculating the Kleene star matrix $\bm{B}^{\ast}$ with a time complexity of at most $O(n^{4})$, when computed by direct matrix multiplications, which can be reduced to $O(n^{3})$ by using the Floyd–Warshall algorithm. Given the matrix $\bm{B}^{\ast}$, each of three terms in the big brackets on the right-hand side of \eqref{E-theta} takes time of $O(n^{2})$, and thus the overall time to compute $\theta$ is no more than $O(m^{2}n^{2})$. 

Note that problem \eqref{P-minhjrjxxrjwj-rjxxrjdj-Bxx-fxg} described in the $(\max,+)$-algebra setting can be solved as a linear program using a polynomial-time iterative procedure such as the Karmarkar algorithm. However, this approach can offer a numerical solution rather than a complete, direct solution in an analytical form like that provided by Theorem~\ref{T-minhjrjxxrjwj-rjxxrjdj-Bxx-fxg}.  

Finally, we represent the result of Theorem~\ref{T-minhjrjxxrjwj-rjxxrjdj-Bxx-fxg} in terms of conventional algebra. For the matrix $\bm{B}$, we denote the entries of the matrix $\bm{B}^{\ast}$ as $b_{ik}^{\ast}$ and note that
\begin{equation*}
b_{ik}^{\ast}
=
\begin{cases}
\beta_{ik},
&
\text{if $i\ne k$},
\\
\max\{\beta_{ik},0\},
&
\text{if $i=k$};
\end{cases}
\end{equation*}
where
\begin{equation*}
\beta_{ik}
=
\max_{1\leq l\leq n-1}\max_{\substack{1\leq i_{1},\ldots,i_{l-1}\leq n\\i_{0}=i,\ i_{l}=k}}(b_{i_{0}i_{1}}+\cdots+b_{i_{l-1}i_{l}}),
\qquad
i,k=1,\ldots,n.
\end{equation*}

With the identity $\max(a,b)=-\min(-a,-b)$ used to save writing, we arrive at the next corollary. 
\begin{corollary}
Suppose that the following conditions hold:
\begin{equation*}
\begin{split}
\text{\upshape 1.}
&\
\max_{\substack{1\leq i_{1},\ldots,i_{k-1}\leq n\\i_{0}=i_{k}=i}}(b_{i_{0}i_{1}}+\cdots+b_{i_{k-1}i_{k}})\leq0,
\qquad
i,k=1,\ldots,n;
\\
\text{\upshape 2.}
&\
b_{ik}^{\ast}
+
\max\left\{
\max_{1\leq l\leq m}(p_{kl}-d_{l}),f_{k}
\right\}
\leq
\min\left\{
\min_{1\leq j\leq m}(p_{ij}+d_{j}),g_{i}
\right\},
\qquad
i,k=1,\ldots,n.
\end{split}
\end{equation*}

Then, the minimum value in problem \eqref{P-minhjrjxxrjwj-rjxxrjdj-Bxx-fxg} is equal to
\begin{multline*}
\theta
=
\max_{1\leq i,k\leq n}
\max_{1\leq j,l\leq m}
\max
\Bigg\{
\frac{w_{l}h_{j}}{w_{j}+w_{l}}+\frac{w_{j}h_{l}}{w_{j}+w_{l}}+\frac{w_{j}w_{l}}{w_{j}+w_{l}}(b_{ik}^{\ast}-p_{ij}+p_{kl}),
\\
h_{j}+w_{j}(b_{ik}^{\ast}-p_{ij}+\max\{p_{kl}-d_{l},f_{k}\}),
\\
h_{l}+w_{l}(b_{ik}^{\ast}-\min\{p_{ij}+d_{i},g_{i}\}+p_{kl})
\Bigg\}.
\end{multline*}

All solution vectors $\bm{x}=(x_{i})$ are given in parametric form by
\begin{equation*}
x_{i}
=
\max_{1\leq k\leq n}(b_{ik}^{\ast}+u_{k}),
\qquad
i=1,\ldots,n;
\end{equation*}
where the vector of parameters $\bm{u}=(u_{k})$ satisfies the condition
\begin{multline*}
\max_{1\leq j\leq m}
\max
\left\{
\frac{h_{j}-\theta}{w_{j}}+p_{kj},
p_{kj}-d_{j},
f_{k}
\right\}
\leq
u_{k}
\\
\leq
\min_{1\leq i\leq n}
\min_{1\leq j\leq m}
\left(
\min
\left\{
\frac{\theta-h_{j}}{w_{j}}+p_{ij},
p_{ij}+d_{j},
g_{i}
\right\}
-
b_{ik}^{\ast}
\right),
\qquad
k=1,\ldots,n.
\end{multline*}
\end{corollary}

To conclude, we note that two conditions of the corollary are responsible for the consistency of the set of inequality constraints that define a set of half-spaces and a hyper-rectangle. The first condition assumes the half-spaces to have a nontrivial intersection, whereas the second implies that this intersection has common points with the hyper-rectangle.

\subsection{Solution for the General Case}

We now suppose that the feasible location set is defined in general form at \eqref{E-S-bikxkckxici_1} and consider the problem
\begin{equation}
\begin{aligned}
&
\text{min}
&&
\bigoplus_{1\leq j\leq m}h_{j}(\bm{p}_{j}^{-}\bm{x}\oplus\bm{x}^{-}\bm{p}_{j})^{w_{j}};
\\
&
\text{s.~t.}
&&
\bm{p}_{j}^{-}\bm{x}
\oplus
\bm{x}^{-}\bm{p}_{j}
\leq
d_{j},
\quad
j=1,\ldots,m;
\\
&&&
b_{ik}x_{k}^{c_{k}}
\leq
x_{i}^{c_{i}},
\quad
i,k=1,\ldots,n;
\\
&&&
\bm{f}
\leq
\bm{x}
\leq\bm{g}.
\end{aligned}
\label{P-minhjrjxxrjwj-rjxxrjdj-bikxkckxici-fxg}
\end{equation}

As in the previous case, we introduce a parameter $\theta$ that represents the minimum value of the objective function in the problem. We reduce the problem to a parametrized system of inequalities, which, after transformation to double inequalities, takes the form
\begin{equation*}
\begin{aligned}
\theta^{-1/w_{j}}h_{j}^{1/w_{j}}\bm{p}_{j}
&\leq
\bm{x}
\leq
\theta^{1/w_{j}}h_{j}^{-1/w_{j}}\bm{p}_{j},
\quad
j=1,\ldots,m;
\\
d_{j}^{-1}\bm{p}_{j}
&\leq
\bm{x}
\leq
d_{j}\bm{p}_{j},
\quad
j=1,\ldots,m;
\\
b_{ik}x_{k}^{c_{k}}
&\leq
x_{i}^{c_{i}},
\quad
i,k=1,\ldots,n;
\\
\bm{f}
&\leq
\bm{x}
\leq
\bm{g}.
\end{aligned}
\end{equation*}

To handle this system, we first rewrite it in the scalar form 
\begin{equation*}
\begin{aligned}
\theta^{-1/w_{j}}h_{j}^{1/w_{j}}p_{ij}
&\leq
x_{i}
\leq
\theta^{1/w_{j}}h_{j}^{-1/w_{j}}p_{ij},
\quad
i=1,\ldots,n;
\quad
j=1,\ldots,m;
\\
d_{j}^{-1}p_{ij}
&\leq
x_{i}
\leq
d_{j}p_{ij},
\quad
i=1,\ldots,n;
\quad
j=1,\ldots,m;
\\
b_{ik}x_{k}^{c_{k}}
&\leq
x_{i}^{c_{i}},
\quad
i,k=1,\ldots,n;
\\
f_{i}
&\leq
x_{i}
\leq
g_{i},
\quad
i=1,\ldots,n.
\end{aligned}
\end{equation*}

Next, we introduce new variables
\begin{equation*}
y_{i}
=
x_{i}^{c_{i}},
\quad
i=1,\ldots,n;
\end{equation*}
and note that $x_{i}=y_{i}^{1/c_{i}}$ for all $i=1,\ldots,n$.

Changing variables for all $i$ such that $c_{i}>0$ yields the system of inequalities
\begin{equation*}
\begin{aligned}
\theta^{-c_{i}/w_{j}}h_{j}^{c_{i}/w_{j}}p_{ij}^{c_{i}}
&\leq
y_{i}
\leq
\theta^{c_{i}/w_{j}}h_{j}^{-c_{i}/w_{j}}p_{ij}^{c_{i}},
\quad
i=1,\ldots,n;
\quad
j=1,\ldots,m;
\\
d_{j}^{-c_{i}}p_{ij}^{c_{i}}
&\leq
y_{i}
\leq
d_{j}^{c_{i}}p_{ij}^{c_{i}},
\quad
i=1,\ldots,n;
\quad
j=1,\ldots,m;
\\
b_{ik}y_{k}
&\leq
y_{i},
\quad
i,k=1,\ldots,n;
\\
f_{i}^{c_{i}}
&\leq
y_{i}
\leq
g_{i}^{c_{i}},
\quad
i=1,\ldots,n;
\end{aligned}
\end{equation*}
whereas for all $i$ with $c_{i}<0$ yields
\begin{equation*}
\begin{aligned}
\theta^{c_{i}/w_{j}}h_{j}^{-c_{i}/w_{j}}p_{ij}^{c_{i}}
&\leq
y_{i}
\leq
\theta^{-c_{i}/w_{j}}h_{j}^{c_{i}/w_{j}}p_{ij}^{c_{i}},
\quad
i=1,\ldots,n;
\quad
j=1,\ldots,m;
\\
d_{j}^{c_{i}}p_{ij}^{c_{i}}
&\leq
y_{i}
\leq
d_{j}^{-c_{i}}p_{ij}^{c_{i}},
\quad
i=1,\ldots,n;
\quad
j=1,\ldots,m;
\\
b_{ik}y_{k}
&\leq
y_{i},
\quad
i,k=1,\ldots,n;
\\
g_{i}^{c_{i}}
&\leq
y_{i}
\leq
f_{i}^{c_{i}},
\quad
i=1,\ldots,n.
\end{aligned}
\end{equation*}

We can summarize both systems as follows:
\begin{equation*}
\begin{aligned}
\theta^{-|c_{i}|/w_{j}}h_{j}^{|c_{i}|/w_{j}}p_{ij}^{c_{i}}
&\leq
y_{i}
\leq
\theta^{|c_{i}|/w_{j}}h_{j}^{-|c_{i}|/w_{j}}p_{ij}^{c_{i}},
\quad
i=1,\ldots,n;
\quad
j=1,\ldots,m;
\\
d_{j}^{-|c_{i}|}p_{ij}^{c_{i}}
&\leq
y_{i}
\leq
d_{j}^{|c_{i}|}p_{ij}^{c_{i}},
\quad
i=1,\ldots,n;
\quad
j=1,\ldots,m;
\\
b_{ik}y_{k}
&\leq
y_{i},
\quad
i,k=1,\ldots,n;
\\
(f_{i}^{-c_{i}}\oplus g_{i}^{-c_{i}})^{-1}
&\leq
y_{i}
\leq
f_{i}^{c_{i}}\oplus g_{i}^{c_{i}},
\quad
i=1,\ldots,n.
\end{aligned}
\end{equation*}

Combining the inequalities into one, we obtain the system of double inequalities 
\begin{multline*}
b_{ik}y_{k}
\oplus
\bigoplus_{1\leq j\leq m}
\theta^{-|c_{i}|/w_{j}}h_{j}^{|c_{i}|/w_{j}}p_{ij}^{c_{i}}
\oplus
\bigoplus_{1\leq j\leq m}
d_{j}^{-|c_{i}|}p_{ij}^{c_{i}}
\oplus
(f_{i}^{-c_{i}}\oplus g_{i}^{-c_{i}})^{-1}
\leq
y_{i}
\\
\leq
\left(
\bigoplus_{1\leq j\leq m}
\theta^{-|c_{i}|/w_{j}}h_{j}^{|c_{i}|/w_{j}}p_{ij}^{-c_{i}}
\oplus
\bigoplus_{1\leq j\leq m}
d_{j}^{-|c_{i}|}p_{ij}^{-c_{i}}
\oplus
(f_{i}^{c_{i}}\oplus g_{i}^{c_{i}})^{-1}
\right)^{-1},
\quad
i,k=1,\ldots,n.
\end{multline*}

We introduce the matrix and vector notation
\begin{equation*}
\bm{B}
=
(b_{ik}),
\qquad
\bm{y}
=
(y_{i}),
\qquad
\bm{q}
=
(q_{i}),
\qquad
\bm{r}
=
(r_{i}),
\qquad
\bm{s}
=
(s_{i}),
\qquad
\bm{t}
=
(t_{i}),
\end{equation*}
with the entries of the last four vectors given by
\begin{equation}
\begin{aligned}
q_{i}
&=
\bigoplus_{1\leq j\leq m}
\theta^{-|c_{i}|/w_{j}}
h_{j}^{|c_{i}|/w_{j}}
p_{ij}^{c_{i}},
&
r_{i}^{-1}
&=
\bigoplus_{1\leq j\leq m}
\theta^{-|c_{i}|/w_{j}}
h_{j}^{|c_{i}|/w_{j}}
p_{ij}^{-c_{i}},
\\
s_{i}
&=
\bigoplus_{1\leq j\leq m}
d_{j}^{-|c_{i}|}
p_{ij}^{c_{i}}
\oplus
(f_{i}^{-c_{i}}\oplus g_{i}^{-c_{i}})^{-1},
&
t_{i}^{-1}
&=
\bigoplus_{1\leq j\leq m}
d_{j}^{-|c_{i}|}
p_{ij}^{-c_{i}}
\oplus
(f_{i}^{c_{i}}\oplus g_{i}^{c_{i}})^{-1},
\end{aligned}
\label{E-qc-rc-sc-tc}
\end{equation}
where $\theta$ represents the minimum value of the objective function.

With this notation, the system of double inequalities can be written as the vector inequality
\begin{equation*}
\bm{B}\bm{y}
\oplus
(\bm{q}\oplus\bm{s})
\leq
\bm{y}
\leq
(\bm{r}^{-}\oplus\bm{t}^{-})^{-},
\end{equation*}
which takes the form of \eqref{I-Bxps-x-xqt} with $\bm{y}$ in place of $\bm{x}$.

Using the same arguments as before, under the condition that $\mathop\mathrm{Tr}(\bm{B})\leq\mathbb{1}$, we represent the solution of the inequality in the parametric form
\begin{equation*}
\bm{y}
=
\bm{B}^{\ast}\bm{v},
\qquad
\bm{q}\oplus\bm{s}
\leq
\bm{v}
\leq
((\bm{r}^{-}\oplus\bm{t}^{-})\bm{B}^{\ast})^{-}.
\end{equation*}

To evaluate the minimum $\theta$ of the objective function and to obtain additional conditions for consistency of constraints, we need to examine the inequality 
\begin{equation*}
\bm{q}\oplus\bm{s}
\leq
((\bm{r}^{-}\oplus\bm{t}^{-})\bm{B}^{\ast})^{-}.
\end{equation*}

After rearrangement of terms, we represent the inequality as the system
\begin{equation*}
\bm{r}^{-}\bm{B}^{\ast}\bm{q}
\leq
\mathbb{1},
\qquad
\bm{r}^{-}\bm{B}^{\ast}\bm{s}
\leq
\mathbb{1},
\qquad
\bm{t}^{-}\bm{B}^{\ast}\bm{q}
\leq
\mathbb{1},
\qquad
\bm{t}^{-}\bm{B}^{\ast}\bm{s}
\leq
\mathbb{1}
\end{equation*}
of the same form as in the previous case, but with the vectors given by \eqref{E-qc-rc-sc-tc}. 

Next, we expand the first three inequalities and write them in scalar form to have
\begin{align*}
(
\theta^{-|c_{i}|/w_{j}}
h_{j}^{|c_{i}|/w_{j}}
p_{ij}^{-c_{i}}
)
b_{ik}^{\ast}
(
\theta^{-|c_{k}|/w_{l}}
h_{l}^{|c_{k}|/w_{l}}
p_{kl}^{c_{k}}
)
&\leq
\mathbb{1},
\quad
i,k=1,\ldots,n;
\quad
j,l=1,\ldots,m;
\\
(
\theta^{-|c_{i}|/w_{j}}
h_{j}^{|c_{i}|/w_{j}}
p_{ij}^{-c_{i}}
)
b_{ik}^{\ast}
s_{k}
&\leq
\mathbb{1},
\quad
i,k=1,\ldots,n;
\quad
j=1,\ldots,m;
\\
t_{i}^{-1}
b_{ik}^{\ast}
(
\theta^{-|c_{k}|/w_{l}}
h_{l}^{|c_{k}|/w_{l}}
p_{kl}^{c_{k}}
)
&\leq
\mathbb{1},
\quad
i,k=1,\ldots,n;
\quad
l=1,\ldots,m;
\\
\bm{t}^{-}\bm{B}^{\ast}\bm{s}
&\leq
\mathbb{1}.
\end{align*}

We solve the first three inequalities with respect to $\theta$ to obtain
\begin{align*}
(
h_{j}^{|c_{i}|/w_{j}}
p_{ij}^{-c_{i}}
b_{ik}^{\ast}
h_{l}^{|c_{k}|/w_{l}}
p_{kl}^{c_{k}}
)^{\frac{w_{j}w_{l}}{|c_{i}|w_{l}+|c_{k}|w_{j}}}
&\leq
\theta,
\quad
i,k=1,\ldots,n;
\quad
j,l=1,\ldots,m;
\\
h_{j}
(
p_{ij}^{-c_{i}}
b_{ik}^{\ast}
s_{k}
)^{w_{j}/|c_{i}|}
&\leq
\theta,
\quad
i,k=1,\ldots,n;
\quad
j=1,\ldots,m;
\\
h_{l}
(
t_{i}^{-1}
b_{ik}^{\ast}
p_{kl}^{c_{k}}
)^{w_{l}/|c_{k}|}
&\leq
\theta,
\quad
i,k=1,\ldots,n;
\quad
l=1,\ldots,m;
\\
\bm{t}^{-}\bm{B}^{\ast}\bm{s}
&\leq
\mathbb{1},
\end{align*}
and then combine these inequalities to write
\begin{gather*}
\theta
\geq
\bigoplus_{1\leq i,k\leq n}
\bigoplus_{1\leq j,l\leq m}
\bigg(
h_{j}^{\frac{|c_{i}|w_{l}}{|c_{i}|w_{l}+|c_{k}|w_{j}}}
h_{l}^{\frac{|c_{k}|w_{j}}{|c_{i}|w_{l}+|c_{k}|w_{j}}}
(
p_{ij}^{-c_{i}}
b_{ik}^{\ast}
p_{kl}^{c_{k}}
)^{\frac{w_{j}w_{l}}{|c_{i}|w_{l}+|c_{k}|w_{j}}}
\\
\oplus
h_{j}
(
p_{ij}^{-c_{i}}
b_{ik}^{\ast}
s_{k}
)^{w_{j}/|c_{i}|}
\oplus
h_{l}
(
t_{i}^{-1}
b_{ik}^{\ast}
p_{kl}^{c_{k}}
)^{w_{l}/|c_{k}|}
\bigg),
\\
\bm{t}^{-}\bm{B}^{\ast}\bm{s}
\leq
\mathbb{1}.
\end{gather*}

The first inequality rewritten as equality yields the minimum value of the objective function, whereas the second inequality gives an additional condition for the consistency of the constraints. 

The result obtained can be formulated as the following statement.

\begin{theorem}
\label{T-minhjrjxxrjwj-rjxxrjdj-bikxkckxici-fxg}
With notation \eqref{E-qc-rc-sc-tc}, suppose that the following conditions hold:
\begin{enumerate}
\item
$\mathop\mathrm{Tr}(\bm{B})\leq\mathbb{1}$,
\item
$\bm{t}^{-}\bm{B}^{\ast}\bm{s}\leq\mathbb{1}$.
\end{enumerate}

Then, the minimum value in problem \eqref{P-minhjrjxxrjwj-rjxxrjdj-bikxkckxici-fxg} is equal to
\begin{multline*}
\theta
=
\bigoplus_{1\leq i,k\leq n}
\bigoplus_{1\leq j,l\leq m}
\bigg(
h_{j}^{\frac{|c_{i}|w_{l}}{|c_{i}|w_{l}+|c_{k}|w_{j}}}
h_{l}^{\frac{|c_{k}|w_{j}}{|c_{i}|w_{l}+|c_{k}|w_{j}}}
(
p_{ij}^{-c_{i}}
b_{ik}^{\ast}
p_{kl}^{c_{k}}
)^{\frac{w_{j}w_{l}}{|c_{i}|w_{l}+|c_{k}|w_{j}}}
\\
\oplus
h_{j}
(
p_{ij}^{-c_{i}}
b_{ik}^{\ast}
s_{k}
)^{w_{j}/|c_{i}|}
\oplus
h_{l}
(
t_{i}^{-1}
b_{ik}^{\ast}
p_{kl}^{c_{k}}
)^{w_{l}/|c_{k}|}
\bigg).
\end{multline*}

All solution vectors $\bm{x}=(x_{i})$ have the elements
\begin{equation*}
x_{i}
=
y_{i}^{1/c_{i}},
\qquad
i=1,\ldots,n,
\end{equation*}
defined by the elements of the vector $\bm{y}=(y_{i})$ which is given by
\begin{equation*}
\bm{y}
=
\bm{B}^{\ast}\bm{v},
\end{equation*}
where the vector of parameters $\bm{v}$ satisfies the condition
\begin{equation*} 
\bm{q}\oplus\bm{s}
\leq
\bm{v}
\leq
((\bm{r}^{-}\oplus\bm{t}^{-})\bm{B}^{\ast})^{-}.
\end{equation*}
\end{theorem}

After translating back into the language of conventional algebra, the result takes the following form.
\begin{corollary}
Suppose that the following conditions hold:
\begin{equation*}
\begin{split}
\text{\upshape 1.}
&\
\max_{\substack{1\leq i_{1},\ldots,i_{k-1}\leq n\\i_{0}=i_{k}=i}}(b_{i_{0}i_{1}}+\cdots+b_{i_{k-1}i_{k}})\leq0,
\qquad
i,k=1,\ldots,n;
\\
\text{\upshape 2.}
&\
b_{ik}^{\ast}
+
\max\left\{
\max_{1\leq l\leq m}(c_{k}p_{kl}-|c_{k}|d_{l}),
\min\{c_{k}f_{k},c_{k}g_{k}\}
\right\}
\\
&\qquad
\leq
\min\left\{
\min_{1\leq j\leq m}(c_{i}p_{ij}+|c_{i}|d_{j}),
\max\{c_{i}f_{i},c_{i}g_{i}\}
\right\},
\qquad
i,k=1,\ldots,n.
\end{split}
\end{equation*}

Then, the minimum value in problem \eqref{P-minhjrjxxrjwj-rjxxrjdj-bikxkckxici-fxg} is equal to
\begin{multline*}
\theta
=
\max_{1\leq i,k\leq n}
\max_{1\leq j,l\leq m}
\max
\Bigg\{
\frac{|c_{i}|w_{l}h_{j}}{|c_{i}|w_{l}+|c_{k}|w_{j}}
+
\frac{|c_{k}|w_{j}h_{l}}{|c_{i}|w_{l}+|c_{k}|w_{j}}
\\
+
\frac{w_{j}w_{l}}{|c_{k}|w_{j}+|c_{i}|w_{l}}(b_{ik}^{\ast}-c_{i}p_{ij}+c_{k}p_{kl}),
\\
h_{j}+\frac{w_{j}}{|c_{i}|}(b_{ik}^{\ast}-c_{i}p_{ij}+\max\{c_{k}p_{kl}-|c_{k}|d_{l},\min\{c_{k}f_{k},c_{k}g_{k}\}\}),
\\
h_{l}+\frac{w_{l}}{|c_{k}|}(b_{ik}^{\ast}-\min\{c_{i}p_{ij}+|c_{i}|d_{j},\max\{c_{i}f_{i},c_{i}g_{i}\}\}+c_{k}p_{kl})
\Bigg\}.
\end{multline*}

All solution vectors $\bm{x}=(x_{i})$ have the elements
\begin{equation*}
x_{i}
=
y_{i}/c_{i},
\qquad
i=1,\ldots,n,
\end{equation*}
defined by the elements of the vector $\bm{y}=(y_{i})$ which is given by
\begin{equation*}
y_{i}
=
\max_{1\leq k\leq n}(b_{ik}^{\ast}+v_{k}),
\qquad
i=1,\ldots,n,
\end{equation*}
where the vector of parameters $\bm{v}=(v_{k})$ satisfies the condition
\begin{multline*}
\left(
\max_{1\leq j\leq m}
\max
\left\{
\frac{|c_{k}|(h_{j}-\theta)}{w_{j}}+c_{k}p_{kj},
c_{k}p_{kj}-|c_{k}|d_{j},
\min\{c_{k}f_{k},c_{k}g_{k}\}
\right\}
\right)
\leq
v_{k}
\\
\leq
\min_{1\leq i\leq n}
\left(
\min_{1\leq j\leq m}
\min
\left\{
\frac{|c_{i}|(\theta-h_{j})}{w_{j}}+c_{i}p_{ij},
c_{i}p_{ij}+|c_{i}|d_{j},
\max\{c_{i}f_{i},c_{i}g_{i}\}
\right\}
-
b_{ik}^{\ast}
\right),
\\
k=1,\ldots,n.
\end{multline*}
\end{corollary}

Note that, as in the previous case, the conditions of the corollary serve to avoid inconsistent sets of constraints, which make the problem unsolvable.

\section{Location with Rectilinear Distance}

We now turn to the solution of location problems defined on the plane with rectilinear distance, inside rectilinear and tilted strips. The rectilinear distance between two vectors $\bm{r}=(r_{1},r_{2})^{T}$ and $\bm{s}=(s_{1},s_{2})^{T}$ in $\mathbb{R}^{2}$ is given in terms of $(\max,+)$-algebra by
\begin{equation*}
\mathrm{d}_{1}(\bm{r},\bm{s})
=
(s_{1}^{-1}r_{1}\oplus r_{1}^{-1}s_{1})(s_{2}^{-1}r_{2}\oplus r_{2}^{-1}s_{2}).
\end{equation*}

To solve the problems, we extend and further develop the technique proposed in \cite{Krivulin2017Using} to solve unweighted two-dimensional rectilinear location problems. The technique involves the representation of the problem in the form of a tropical optimization problem, followed by a change of variables, which reduces the optimization problem to problems in the form of \eqref{P-minhjrjxxrjwj-rjxxrjdj-Bxx-fxg} and \eqref{P-minhjrjxxrjwj-rjxxrjdj-bikxkckxici-fxg}.

\subsection{Location in Rectilinear Strip}

We represent the feasible location area inside a vertical rectilinear strip, given by \eqref{E-S-f1x2x1g1x2f2x1x2g2x1ax1b}, in the $(\max,+)$-algebra setting as follows
\begin{equation*}
S
=
\{(x_{1},x_{2})^{T}|\ f_{1}x_{2}^{-1}\leq x_{1}\leq g_{1}x_{2}^{-1},\ f_{2}x_{1}\leq x_{2}\leq g_{2}x_{1},\ a\leq x_{1}\leq b\}.
\end{equation*}

After rewriting the distances, location problem \eqref{P-minwjdxrjhj-dxrjdj-xS} takes the form
\begin{equation}
\begin{aligned}
&
\text{min}
&&
\bigoplus_{1\leq j\leq m}h_{j}((p_{1j}^{-1}x_{1}\oplus x_{1}^{-1}p_{1j})(p_{2j}^{-1}x_{2}\oplus x_{2}^{-1}p_{2j}))^{w_{j}};
\\
&
\text{s.~t.}
&&
(p_{1j}^{-1}x_{1}\oplus x_{1}^{-1}p_{1j})(p_{2j}^{-1}x_{2}\oplus x_{2}^{-1}p_{2j})
\leq
d_{j},
\quad
j=1,\ldots,m;
\\
&&&
f_{1}x_{1}^{-1}\leq x_{2}\leq g_{1}x_{1}^{-1},
\quad
f_{2}x_{2}\leq x_{1}\leq g_{2}x_{2},
\quad
a\leq x_{1}\leq b.
\end{aligned}
\label{P-minhjr1jx1x1r1jr2jx2x2r2jwj-r1jx1x1r1jr2jx2x2r2jdj-f1x1x2g1x1-f2x2x1g2x2-ax1b}
\end{equation}

As before, we assume that all parameters and vectors involved in the problem formulation have nonzero values in the sense of $(\max,+)$-algebra. 

The solution of problem \eqref{P-minhjr1jx1x1r1jr2jx2x2r2jwj-r1jx1x1r1jr2jx2x2r2jdj-f1x1x2g1x1-f2x2x1g2x2-ax1b} is based on changing variables to reduce the problem to \eqref{P-minhjrjxxrjwj-rjxxrjdj-Bxx-fxg} and thus to take advantage of the above-obtained results. Note that this transformation from \eqref{P-minhjr1jx1x1r1jr2jx2x2r2jwj-r1jx1x1r1jr2jx2x2r2jdj-f1x1x2g1x1-f2x2x1g2x2-ax1b} to \eqref{P-minhjrjxxrjwj-rjxxrjdj-Bxx-fxg} reflects the known relationship between the location problems on the plane with rectilinear and Chebyshev distances \cite{Dearing1972Onsome,Francis1972Ageometrical}, which can be converted from one to the other by a coordinate transformation.

To solve problem \eqref{P-minhjr1jx1x1r1jr2jx2x2r2jwj-r1jx1x1r1jr2jx2x2r2jdj-f1x1x2g1x1-f2x2x1g2x2-ax1b}, we first introduce new vectors
\begin{equation}
\begin{aligned}
\bm{y}
&=
\left(
\begin{array}{c}
y_{1}
\\
y_{2}
\end{array}
\right),
&
y_{1}
&=
x_{1}x_{2},
&
y_{2}
&=
x_{1}^{-1}x_{2};
\\
\bm{o}_{j}
&=
\left(
\begin{array}{c}
o_{1j}
\\
o_{2j}
\end{array}
\right),
&
o_{1j}
&=
p_{1j}p_{2j},
&
o_{2j}
&=
p_{1j}^{-1}p_{2j},
\quad
j=1,\ldots,m.
\end{aligned}
\label{E-y-oj}
\end{equation}

Clearly, the elements of the vector $\bm{x}$ are related with those of $\bm{y}$ by the equalities
\begin{equation*}
x_{1}
=
y_{1}^{1/2}y_{2}^{-1/2},
\qquad
x_{2}
=
y_{1}^{1/2}y_{2}^{1/2}.
\end{equation*}

With the new notation, we can write
\begin{multline*}
(p_{1j}^{-1}x_{1}\oplus x_{1}^{-1}p_{1j})(p_{2j}^{-1}x_{2}\oplus x_{2}^{-1}p_{2j})
=
o_{1j}^{-1}y_{1}
\oplus
o_{2j}^{-1}y_{2}
\oplus
o_{2j}y_{2}^{-1}
\oplus
o_{1j}y_{1}^{-1}
\\=
\bm{o}_{j}^{-}\bm{y}\oplus\bm{y}^{-}\bm{o}_{j},
\quad
j=1,\ldots,m.
\end{multline*}

It remains to rewrite the constraints which determine the feasible location area $S$. The first inequality $f_{1}x_{2}^{-1}\leq x_{1}\leq g_{1}x_{2}^{-1}$ is equivalent to $f_{1}\leq x_{1}x_{2}\leq g_{1}$ which can be written as $f_{1}\leq y_{1}\leq g_{1}$. In the same way, we represent the second inequality $f_{2}x_{1}\leq x_{2}\leq g_{2}x_{1}$ as $f_{2}\leq y_{2}\leq g_{2}$.

We take the last inequality $a\leq x_{1}\leq b$ and put its left part in the equivalent form $a^{2}x_{1}^{-1}x_{2}\leq x_{1}x_{2}$ which can be expressed as $a^{2}y_{2}\leq y_{1}$. The right part evolves into $x_{1}x_{2}\leq b^{2}x_{1}^{-1}x_{2}$ and then into $b^{-2}y_{1}\leq y_{2}$.

With the vector and matrix notation
\begin{equation}
\bm{f}
=
\left(
\begin{array}{c}
f_{1}
\\
f_{2}
\end{array}
\right),
\qquad
\bm{g}
=
\left(
\begin{array}{c}
g_{1}
\\
g_{2}
\end{array}
\right),
\qquad
\bm{B}
=
\left(
\begin{array}{cc}
\mathbb{0} & a^{2}
\\
b^{-2} & \mathbb{0}
\end{array}
\right),
\label{E-f-g-B}
\end{equation}
where $\mathbb{0}=-\infty$,  we express the constraints in vector form as
\begin{equation*}
\bm{f}
\leq
\bm{y}
\leq
\bm{g},
\qquad
\bm{B}\bm{y}
\leq
\bm{y}.
\end{equation*}

We now obtain problem \eqref{P-minhjr1jx1x1r1jr2jx2x2r2jwj-r1jx1x1r1jr2jx2x2r2jdj-f1x1x2g1x1-f2x2x1g2x2-ax1b} formulated in terms of $(\max,+)$-algebra as follows
\begin{equation}
\begin{aligned}
&
\text{min}
&&
\bigoplus_{1\leq j\leq m}h_{j}(\bm{o}_{j}^{-}\bm{y}\oplus\bm{y}^{-}\bm{o}_{j})^{w_{j}};
\\
&
\text{s.~t.}
&&
\bm{o}_{j}^{-}\bm{y}\oplus\bm{y}^{-}\bm{o}_{j}
\leq
d_{j},
\quad
j=1,\ldots,m;
\\
&&&
\bm{B}\bm{y}
\leq
\bm{y},
\quad
\bm{f}\leq\bm{y}\leq\bm{g}.
\end{aligned}
\label{P-minhjojyyojwj-ojyyojdj-Byy-fyg}
\end{equation}

Since the problem obtained takes the form of \eqref{P-minhjrjxxrjwj-rjxxrjdj-Bxx-fxg}, we apply Theorem~\ref{T-minhjrjxxrjwj-rjxxrjdj-Bxx-fxg} to derive a complete solution. First, we note that, under the condition $a\leq b$, we have
\begin{equation*}
\mathop\mathrm{Tr}(\bm{B})
=
ab^{-1}\leq\mathbb{1},
\qquad
\bm{B}^{\ast}
=
\left(
\begin{array}{cc}
\mathbb{1} & a^{2}
\\
b^{-2} &\mathbb{1}
\end{array}
\right).
\end{equation*}

Taking into account \eqref{E-y-oj}, we rewrite notation \eqref{E-p-q-f-g} as follows:
\begin{equation*}
\begin{aligned}
\bm{q}
&=
\bigoplus_{1\leq j\leq m}
\theta^{-1/w_{j}}h_{j}^{1/w_{j}}
\bm{o}_{j},
&
\bm{r}^{-}
&=
\bigoplus_{1\leq j\leq m}
\theta^{-1/w_{j}}h_{j}^{1/w_{j}}
\bm{o}_{j}^{-},
\\
\bm{s}
&=
\bigoplus_{1\leq j\leq m}
d_{j}^{-1}
\bm{o}_{j}
\oplus
\bm{f},
&
\bm{t}^{-}
&=
\bigoplus_{1\leq j\leq m}
d_{j}^{-1}
\bm{o}_{j}^{-}
\oplus
\bm{g}^{-}.
\end{aligned}
\label{E-p-q-f-g-1}
\end{equation*}

The next statement provides a complete solution of problem \eqref{P-minhjr1jx1x1r1jr2jx2x2r2jwj-r1jx1x1r1jr2jx2x2r2jdj-f1x1x2g1x1-f2x2x1g2x2-ax1b}.

\begin{theorem}\label{T-minhjojyyojwj-ojyyojdj-Byy-fyg}
With the above notation, suppose that the condition $\bm{t}^{-}\bm{B}^{\ast}\bm{s}\leq\mathbb{1}$ holds. 
Then, the minimum value in problem \eqref{P-minhjr1jx1x1r1jr2jx2x2r2jwj-r1jx1x1r1jr2jx2x2r2jdj-f1x1x2g1x1-f2x2x1g2x2-ax1b} is equal to
\begin{equation}
\theta
=
\bigoplus_{1\leq j,l\leq m}
\left(
h_{j}^{\frac{w_{l}}{w_{j}+w_{l}}}
h_{l}^{\frac{w_{j}}{w_{j}+w_{l}}}
(
\bm{o}_{j}^{-}
\bm{B}^{\ast}
\bm{o}_{l}
)^{\frac{w_{j}w_{l}}{w_{j}+w_{l}}}
\oplus
h_{j}
(
\bm{o}_{j}^{-}
\bm{B}^{\ast}
\bm{s}
)^{w_{j}}
\oplus
h_{l}
(
\bm{t}^{-}
\bm{B}^{\ast}
\bm{o}_{l}
)^{w_{l}}
\right).
\label{E-theta3}
\end{equation}

All solution vectors $\bm{x}=(x_{1},x_{2})^{T}$ have the elements
\begin{equation*}
x_{1}
=
y_{1}^{1/2}y_{2}^{-1/2},
\qquad
x_{2}
=
y_{1}^{1/2}y_{2}^{1/2},
\end{equation*}
defined by the elements of the vector $\bm{y}=(y_{1},y_{2})^{T}$ which is given by
\begin{equation*}
\bm{y}
=
\bm{B}^{\ast}\bm{u},
\end{equation*}
where the vector of parameters $\bm{u}$ satisfies the condition
\begin{equation*} 
\bm{q}\oplus\bm{s}
\leq
\bm{u}
\leq
((\bm{r}^{-}\oplus\bm{t}^{-})\bm{B}^{\ast})^{-}.
\end{equation*}
\end{theorem}

To rewrite the result in terms of ordinary arithmetic operations, we first represent the entries of the matrix $\bm{B}^{\ast}$ as
\begin{equation*}
b_{11}^{\ast}=b_{22}^{\ast}=0,
\qquad
b_{12}^{\ast}=2a,
\qquad
b_{21}^{\ast}=-2b.
\end{equation*}

Now, the result of the theorem reads as follows. 
\begin{corollary}
Let $o_{1j}=p_{1j}+p_{2j}$ and $o_{2j}=p_{2j}-p_{1j}$ for all $j=1,\ldots,m$, and suppose that the following condition holds:
\begin{equation*}
b_{ik}^{\ast}
+
\max_{1\leq l\leq m}
\max\{o_{kl}-d_{l},f_{k}\}
\leq
\min_{1\leq j\leq m}
\min\{o_{ij}+d_{j},g_{i}\},
\qquad
i,k=1,2.
\end{equation*}

Then, the minimum value in problem \eqref{P-minhjojyyojwj-ojyyojdj-Byy-fyg} is equal to
\begin{multline*}
\theta
=
\max_{1\leq j,l\leq m}
\max
\Bigg\{
\frac{w_{l}h_{j}}{w_{j}+w_{l}}+\frac{w_{j}h_{l}}{w_{j}+w_{l}}
+
\frac{w_{j}w_{l}}{w_{j}+w_{l}}
\max_{1\leq i,k\leq2}
(b_{ik}^{\ast}-o_{ij}+o_{kl}),
\\
h_{j}+w_{j}
\max_{1\leq i,k\leq2}
(b_{ik}^{\ast}-o_{ij}+\max\{o_{kl}-d_{l},f_{k}\}),
\\
h_{l}+w_{l}
\max_{1\leq i,k\leq2}
(b_{ik}^{\ast}-\min\{d_{j}+o_{ij},g_{i}\}+o_{kl})
\Bigg\}.
\end{multline*}

All solution vectors $\bm{x}=(x_{1},x_{2})^{T}$ have the elements
\begin{equation*}
x_{1}
=
(y_{1}-y_{2})/2,
\qquad
x_{2}
=
(y_{1}+y_{2})/2,
\end{equation*}
defined by the elements of the vector $\bm{y}=(y_{1},y_{2})^{T}$ which is given by
\begin{equation*}
y_{i}
=
\max\{b_{i1}^{\ast}+u_{1},b_{i2}^{\ast}+u_{2}\},
\qquad
i=1,2,
\end{equation*}
where the vector of parameters $\bm{u}=(u_{1},u_{2})^{T}$ satisfies the conditions
\begin{multline*}
\max_{1\leq j\leq m}
\max
\left\{
\frac{h_{j}-\theta}{w_{j}}+o_{kj},
-d_{j}+o_{kj},
f_{k}
\right\}
\leq
u_{k}
\\
\leq
\min_{1\leq i\leq 2}
\min_{1\leq j\leq m}
\left(
\min
\left\{
\frac{\theta-h_{j}}{w_{j}}+o_{ij},
d_{j}+o_{ij},
g_{i}
\right\}
-
b_{ik}^{\ast}
\right),
\qquad
k=1,2.
\end{multline*}
\end{corollary}

\subsection{Location in Tilted Strip}

In this section, we derive a complete solution for the rectilinear location problem in which the feasible location area is given by 
\begin{equation*}
S
=
\{(x_{1},x_{2})^{T}|\ f_{1}x_{2}^{-1}\leq x_{1}\leq g_{1}x_{2}^{-1},\ f_{2}x_{1}\leq x_{2}\leq g_{2}x_{1},\ ax_{2}\leq x_{1}^{c}\leq bx_{2},\ c\ne1\}.
\end{equation*}

To represent the location problem in terms of $(\max,+)$-algebra, we start with the notation given by \eqref{E-y-oj} and \eqref{E-f-g-B} for the vectors $\bm{o}$, $\bm{y}$, $\bm{f}$, $\bm{g}$, and for the matrix $\bm{B}$. Next, we consider the double inequality $ax_{2}\leq x_{1}^{c}\leq bx_{2}$, and introduce the new variables
\begin{equation*}
c_{1}
=
c-1,
\qquad
c_{2}
=
c+1.
\end{equation*}

Multiplying the left part of this inequality by $x_{1}^{-(c+1)/2}x_{2}^{(c-1)/2}$ yields the inequality $a(x_{1}^{-1}x_{2})^{(c+1)/2}\leq (x_{1}x_{2})^{(c-1)/2}$. We change the variables to rewrite the last inequality as $ay_{2}^{c_{2}/2}\leq y_{1}^{c_{1}/2}$ and then as $a^{2}y_{2}^{c_{2}}\leq y_{1}^{c_{1}}$.

In the same way, we represent the right part of the double inequality as $(x_{1}x_{2})^{(c-1)/2}\leq b(x_{1}^{-1}x_{2})^{(c+1)/2}$ and then as $b^{-2}y_{1}^{c_{1}}\leq y_{2}^{c_{2}}$.

The problem under study now takes the form
\begin{equation}
\begin{aligned}
&
\text{min}
&&
\bigoplus_{1\leq j\leq m}h_{j}(\bm{o}_{j}^{-}\bm{y}\oplus\bm{y}^{-}\bm{o}_{j})^{w_{j}};
\\
&
\text{s.~t.}
&&
\bm{o}_{j}^{-}\bm{y}
\oplus
\bm{y}^{-}\bm{o}_{j}
\leq
d_{j},
\quad
j=1,\ldots,m;
\\
&&&
a^{2}y_{2}^{c_{2}}
\leq
y_{1}^{c_{1}},
\quad
b^{-2}y_{1}^{c_{1}}
\leq
y_{2}^{c_{2}},
\\
&&&
\bm{f}
\leq
\bm{y}
\leq\bm{g}.
\end{aligned}
\label{P-minhjojxxojwj-ojxxojdj-ay2c2y1c1-by1c1y2c2-fyg}
\end{equation}

We redefine notation \eqref{E-qc-rc-sc-tc} for the vectors
\begin{equation*}
\bm{q}
=
\left(
\begin{array}{c}
q_{1}
\\
q_{2}
\end{array}
\right),
\qquad
\bm{r}
=
\left(
\begin{array}{c}
r_{1}
\\
r_{2}
\end{array}
\right),
\qquad
\bm{s}
=
\left(
\begin{array}{c}
s_{1}
\\
s_{2}
\end{array}
\right),
\qquad
\bm{t}
=
\left(
\begin{array}{c}
t_{1}
\\
t_{2}
\end{array}
\right),
\end{equation*}
by setting their elements for each $i=1,2$ as follows:
\begin{align*}
q_{i}
&=
\bigoplus_{1\leq j\leq m}
\theta^{-|c_{i}|/w_{j}}
h_{j}^{|c_{i}|/w_{j}}
o_{ij}^{c_{i}},
&
r_{i}^{-1}
&=
\bigoplus_{1\leq j\leq m}
\theta^{-|c_{i}|/w_{j}}
h_{j}^{|c_{i}|/w_{j}}
o_{ij}^{-c_{i}},
\\
s_{i}
&=
\bigoplus_{1\leq j\leq m}
d_{j}^{-|c_{i}|}
o_{ij}^{c_{i}}
\oplus
(f_{i}^{-c_{i}}\oplus g_{i}^{-c_{i}})^{-1},
&
t_{i}^{-1}
&=
\bigoplus_{1\leq j\leq m}
d_{j}^{-|c_{i}|}
o_{ij}^{-c_{i}}
\oplus
(f_{i}^{c_{i}}\oplus g_{i}^{c_{i}})^{-1}.
\end{align*}

Observing that the problem obtained has the form of \eqref{P-minhjrjxxrjwj-rjxxrjdj-bikxkckxici-fxg}, we apply Theorem~\ref{T-minhjrjxxrjwj-rjxxrjdj-bikxkckxici-fxg}, which leads to the following result.
\begin{theorem}
\label{T-minhjojxxojwj-ojxxojdj-ay2c2y1c1-by1c1y2c2-fyg}
With the above notation, suppose that the condition $\bm{t}^{-}\bm{B}^{\ast}\bm{s}\leq\mathbb{1}$ holds. Then, the minimum value in problem \eqref{P-minhjojxxojwj-ojxxojdj-ay2c2y1c1-by1c1y2c2-fyg} is equal to
\begin{multline*}
\theta
=
\bigoplus_{1\leq i,k\leq 2}
\bigoplus_{1\leq j,l\leq m}
\bigg(
h_{j}^{\frac{|c_{i}|w_{l}}{|c_{i}|w_{l}+|c_{k}|w_{j}}}
h_{l}^{\frac{|c_{k}|w_{j}}{|c_{i}|w_{l}+|c_{k}|w_{j}}}
(
o_{ij}^{-c_{i}}
b_{ik}^{\ast}
o_{kl}^{c_{k}}
)^{\frac{w_{j}w_{l}}{|c_{i}|w_{l}+|c_{k}|w_{j}}}
\\
\oplus
h_{j}
(
o_{ij}^{-c_{i}}
b_{ik}^{\ast}
s_{k}
)^{w_{j}/|c_{i}|}
\oplus
h_{l}
(
t_{i}^{-1}
b_{ik}^{\ast}
o_{kl}^{c_{k}}
)^{w_{l}/|c_{k}|}
\bigg).
\end{multline*}

All solution vectors $\bm{x}=(x_{1},x_{2})^{T}$ have the elements
\begin{equation*}
x_{1}
=
(y_{1}y_{2}^{-1})^{1/2c_{1}},
\qquad
x_{2}
=
(y_{1}y_{2})^{1/2c_{2}},
\end{equation*}
defined by the elements of the vector $\bm{y}=(y_{1},y_{2})^{T}$ which is given by
\begin{equation*}
\bm{y}
=
\bm{B}^{\ast}\bm{v},
\end{equation*}
where the vector of parameters $\bm{v}$ satisfies the condition
\begin{equation*} 
\bm{q}\oplus\bm{s}
\leq
\bm{v}
\leq
((\bm{r}^{-}\oplus\bm{t}^{-})\bm{B}^{\ast})^{-}.
\end{equation*}
\end{theorem}

Returning back to conventional algebra yields the next result. 
\begin{corollary}
Let $o_{1j}=p_{1j}+p_{2j}$ and $o_{2j}=p_{2j}-p_{1j}$ for all $j=1,\ldots,m$, and suppose that the following conditions hold:
\begin{multline*}
b_{ik}^{\ast}
+
\max\left\{
\max_{1\leq l\leq m}(c_{k}o_{kl}-|c_{k}|d_{l}),
\min\{c_{k}f_{k},c_{k}g_{k}\}
\right\}
\\
\leq
\min\left\{
\min_{1\leq j\leq m}(c_{i}o_{ij}+|c_{i}|d_{j}),
\max\{c_{i}f_{i},c_{i}g_{i}\}
\right\},
\qquad
i,k=1,2.
\end{multline*}

Then, the minimum value in problem \eqref{P-minhjojxxojwj-ojxxojdj-ay2c2y1c1-by1c1y2c2-fyg} is equal to
\begin{multline*}
\theta
=
\max_{1\leq i,k\leq 2}
\max_{1\leq j,l\leq m}
\max
\Bigg\{
\frac{|c_{i}|w_{l}h_{j}}{|c_{i}|w_{l}+|c_{k}|w_{j}}
+
\frac{|c_{k}|w_{j}h_{l}}{|c_{i}|w_{l}+|c_{k}|w_{j}}
\\
+
\frac{w_{j}w_{l}}{|c_{k}|w_{j}+|c_{i}|w_{l}}(b_{ik}^{\ast}-c_{i}o_{ij}+c_{k}o_{kl}),
\\
h_{j}+\frac{w_{j}}{|c_{i}|}(b_{ik}^{\ast}-c_{i}o_{ij}+\max\{c_{k}o_{kl}-|c_{k}|d_{l},\min\{c_{k}f_{k},c_{k}g_{k}\}\}),
\\
h_{l}+\frac{w_{l}}{|c_{k}|}(b_{ik}^{\ast}-\min\{c_{i}o_{ij}+|c_{i}|d_{j},\max\{c_{i}f_{i},c_{i}g_{i}\}\}+c_{k}o_{kl})
\Bigg\}.
\end{multline*}

All solution vectors $\bm{x}=(x_{1},x_{2})^{T}$ have the elements
\begin{equation*}
x_{1}
=
(y_{1}-y_{2})/2c_{1},
\qquad
x_{2}
=
(y_{1}+y_{2})/2c_{2},
\end{equation*}
defined by the elements of the vector $\bm{y}=(y_{1},y_{2})^{T}$ which is given by
\begin{equation*}
y_{i}
=
\max\{b_{i1}^{\ast}+v_{1},b_{i2}^{\ast}+v_{2}\},
\qquad
i=1,2,
\end{equation*}
where the vector of parameters $\bm{v}=(v_{1},v_{2})^{T}$ satisfies the conditions
\begin{multline*}
\left(
\max_{1\leq j\leq m}
\max
\left\{
\frac{|c_{k}|(h_{j}-\theta)}{w_{j}}+c_{k}o_{kj},
c_{k}o_{kj}-|c_{k}|d_{j},
\min\{c_{k}f_{k},c_{k}g_{k}\}
\right\}
\right)
\leq
v_{k}
\\
\leq
\min_{1\leq i\leq 2}
\left(
\min_{1\leq j\leq m}
\min
\left\{
\frac{|c_{i}|(\theta-h_{j})}{w_{j}}+c_{i}o_{ij},
c_{i}o_{ij}+|c_{i}|d_{j},
\max\{c_{i}f_{i},c_{i}g_{i}\}
\right\}
-
b_{ik}^{\ast}
\right),
\\
k=1,2.
\end{multline*}
\end{corollary}

\section{Conclusions}

The paper has examined minimax single-facility location problems in the $n$-dimensional vector space with Chebyshev distance and in the two-dimensional plane with rectilinear distance. The feasible  location areas are given by sets of inequality constraints, which define the intersection of half-spaces and a hyper-rectangle for Chebyshev location, and the intersection of half-spaces and a strip area for rectilinear location.

We have started with location problems in a space of arbitrary dimension with Chebyshev distance. To handle the problems, we first represented them in terms of $(\max,+)$-algebra as a tropical optimization problem. The solution approach was implemented, which introduces an additional parameter to represent the optimal value of the objective function, and then uses properties of the operations in $(\max,+)$-algebra to reduce the optimization problem to the solution of a set of parametrized inequalities. The existence conditions for solutions of the system serve to evaluate the parameter, whereas all solutions of the system are taken as a complete solution of the optimization problem.

Using this approach, we have derived new exact, complete solutions of the multidimensional location problems with Chebyshev distance in terms of tropical mathematics, and represented the solutions in the conventional form. The results obtained were extended to examine two-dimensional problems with rectilinear distance and to provide new solutions in both tropical and conventional algebra settings. The solutions are given in a closed form suitable for further analytical study and direct computations with low polynomial complexity in terms of both the dimension of the location space and the number of given points. 

Possible lines of further research include the development of algebraic methods to solve rectilinear location problems in the three-dimensional space and in the space of arbitrary dimension, as well as to solve Chebyshev and rectilinear problems with new types of constraints.

\bibliographystyle{abbrvurl}

\bibliography{Algebraic_solution_of_minimax_single-facility_constrained_location_problems_with_Chebyshev_and_rectilinear_distances}

\end{document}